\pdfoutput=1
\RequirePackage{ifpdf}
\ifpdf 
\documentclass[pdftex]{sigma}
\else
\documentclass{sigma}
\fi

\newcommand{\CC}{\mathbb{C}}
\newcommand{\QQ}{\mathbb{Q}}
\newcommand{\ZZ}{\mathbb Z}

\newcommand{\PP}{\mathbb{P}}
\newcommand{\tv}{\mathbb{TV}}

\newcommand{\D}{\mathcal D}
\newcommand{\cX}{\mathcal{X}}

\DeclareMathOperator{\Proj}{Proj}
\DeclareMathOperator{\spec}{Spec}
\DeclareMathOperator{\conv}{conv}
\DeclareMathOperator{\cone}{cone}
\DeclareMathOperator{\tail}{tail}
\DeclareMathOperator{\Hom}{Hom}
\DeclareMathOperator{\aut}{Aut}

\newtheorem{Theorem}{Theorem}[section]
 { \theoremstyle{definition}
\newtheorem{Definition}[Theorem]{Definition}
\newtheorem{Example}[Theorem]{Example}
\newtheorem{Remark}[Theorem]{Remark} }

\begin{document}

\allowdisplaybreaks

\renewcommand{\thefootnote}{$\star$}

\renewcommand{\PaperNumber}{047}

\FirstPageHeading

\ShortArticleName{Mutations of Laurent Polynomials and Flat Families with Toric Fibers}

\ArticleName{Mutations of Laurent Polynomials\\ and Flat Families with Toric Fibers\footnote{This
paper is a contribution to the Special Issue ``Mirror Symmetry and Related Topics''. The full collection is available at \href{http://www.emis.de/journals/SIGMA/mirror_symmetry.html}{http://www.emis.de/journals/SIGMA/mirror\_{}symmetry.html}}}

\Author{Nathan Owen ILTEN}

\AuthorNameForHeading{N.O.~Ilten}

\Address{Department of Mathematics, University of California, Berkeley CA 94720, USA}
\Email{\href{mailto:nilten@math.berkeley.edu}{nilten@math.berkeley.edu}}
\URLaddress{\url{http://math.berkeley.edu/~nilten/}}

\ArticleDates{Received May 21, 2012, in f\/inal form July 25, 2012; Published online July 28, 2012}

\Abstract{We give a general criterion for two toric varieties to appear as f\/ibers in a f\/lat family over $\PP^1$.
We apply this to show that certain birational transformations mapping a~Laurent polynomial to another Laurent polynomial correspond to deformations between the associated toric varieties.}

\Keywords{toric varieties; mirror symmetry; deformations; Newton polyhedra}
\Classification{14M25; 14D06; 53D37}

\section{Introduction}
Consider the $n$-dimensional torus $T=(\CC^*)^n$ with coordinates $z_1,\ldots,z_n$ along with the logarithmic volume form
\begin{gather*}
\omega=\frac{1}{(2\pi i)^n}\frac{d z_1}{z_1}\wedge \cdots \wedge \frac{d z_n}{z_n},
\end{gather*}
where $z_i$ are coordinates on $T$.
\begin{Definition}[{cf.~\cite[Def\/inition~7]{galkin}}]
	Let $f\colon T\to\CC$ be a Laurent polynomial in $z_1,\ldots,z_n$.
	A~\emph{mutation} of $f$ is a birational transformation $\phi\in\aut(\CC(z_1,\ldots,z_n))$
preserving $\omega$ such that $\phi(f)$ is again a Laurent polynomial.
\end{Definition}

 Such  transformations arise in the context of wall-crossing for counts of holomorphic discs bounded by special Lagrangian tori, see e.g.~\cite{auroux}.

\begin{Example}\label{ex:easymut}
Let $g$ be a non-zero Laurent polynomial in $z_2,\ldots,z_n$, and consider the birational transformation
\begin{gather*}	
	\phi: \ (z_1,z_2,\ldots,z_n)\mapsto (z_1/g,z_2,\ldots,z_n).
\end{gather*}
Then $\phi$ is a mutation of a Laurent polynomial $f$ if and only if $f$ can be written as
\begin{gather*}
f=\sum_{i=k}^l f_i z_1^i
\end{gather*}
with $f_i$ Laurent polynomials in $z_2,\ldots,z_n$ such that, for $i>0$, $f_i/g^i$ is a Laurent polynomial.
\end{Example}

The purpose of this article is to relate these special mutations to deformations. Given a~lattice polytope $\Delta$ containing the origin in its interior, let
$\Sigma(\Delta)$ be the face fan of $\Delta$, that is, the fan whose cones are generated by proper faces of $\Delta$. Let $\tv(\Delta)$ denote the toric variety corresponding to the fan $\Sigma(\Delta)$, see~\cite{fulton:93a} for details on toric varieties.  Our main result is

\begin{Theorem}\label{thm:main}
	Let $\phi$ be a mutation of a Laurent polynomial $f$ of the type of Example~{\rm \ref{ex:easymut}} and suppose that $\Delta(f)$ contains the origin in its interior. Then there is a flat projective family $\pi\colon \cX\to\PP^1$ such that $\pi^{-1}(0)=\tv(\Delta(f))$ and $\pi^{-1}(\infty)=\tv(\Delta(\phi(f))$.
\end{Theorem}

A special case of the mutations from Example \ref{ex:easymut} are those considered by Galkin and Usnich in~\cite{galkin}. In particular, the above theorem shows that any two toric surfaces related via one of the mutations of~\cite{galkin} appear as special f\/ibers in a f\/lat projective family over $\PP^1$, answering a question posed in loc.~cit.

The family $\pi\colon \cX\to\PP^1$  comes from a general construction of R.~Vollmert and the author, see~\cite{ilten:09b}. In Section~\ref{sec:fam}, we use this construction to formulate a criterion for two toric varieties to appear as f\/ibers in a f\/lat family over $\PP^1$, see Theorem~\ref{thm:family}.
We apply this to prove Theorem~\ref{thm:main} in Section~\ref{sec:main}. In Section~\ref{sec:surfaces} we brief\/ly discuss the mutations of Usnich and Galkin.

We end this section by f\/ixing notation. Let $N$ be a lattice, and $M=\Hom(N,\ZZ)$ its dual. We denote the $\QQ$-vector spaces $N\otimes \QQ$ and $M\otimes \QQ$ by $N_\QQ$ and $M_\QQ$, respectively. To any polyhedral cone $\sigma\subset N_\QQ$, we can associate the af\/f\/ine toric variety $\tv(\sigma)=\spec \CC[\sigma^\vee\cap M]$. Likewise, for any polytope $\Delta\subset N_\QQ$ containing the origin in its interior, let $\tv(\Delta)$ denote the toric variety corresponding to the face fan $\Sigma(\Delta)$ of $\Delta$.

Let $\Delta_0$, $\Delta_1$ be polyhedra in $N_\QQ$. Their \emph{Minkowski sum} $\Delta_0+\Delta_1$ consists of all points $v_0+v_1$ where $v_0\in\Delta_0$, $v_1\in\Delta_1$. Given a polyhedron $\Delta\subset N_\QQ$, its \emph{tailcone} $\tail(\Delta)$ consists of all $v\in N_\QQ$ such that $v+\Delta\subset \Delta$.

\begin{Definition}
A pair of polyhedra $\Delta_0$, $\Delta_1$ in $N_\QQ$ is \emph{admissible} if their tailcones are equal, and if for all $u\in M\cap\tail(\Delta_0)^\vee$, the  minimum value of $u$ on $\Delta_i$ is integral for either $i=0$ or $i=1$ (or both).
\end{Definition}

\section{Families with toric f\/ibers}\label{sec:fam}

We now describe a criterion for two toric varieties $\cX_0$, $\cX_\infty$ to appear as f\/ibers in
a f\/lat family $\pi:\cX\to\PP^1$. For now, we will focus on the case where $\cX_0$ and $\cX_\infty$ are af\/f\/ine, as well as making a number of simplifying assumptions.

The idea is to start with $\cX_0$, and provide a procedure for constructing possible $\cX_\infty$. So let us f\/ix some af\/f\/ine toric variety $\cX_0=\tv(\sigma)$, where $\sigma$ is a cone in $N_\QQ$. To construct $\cX_\infty$, we will make three choices, only two of which will have an ef\/fect on the end result. First of all, we will choose a primitive $u\in M=\Hom(N,\ZZ)$ such that $\pm u\notin\sigma^\vee$. This gives rise to an exact sequence
\begin{gather*}
\begin{CD}
0 @>>> N' @>\iota >> N @>u>> \ZZ @>>> 0.
\end{CD}
\end{gather*}
Note that the choice of $u$ corresponds to a subtorus $T'=\CC^*\otimes N'$ of the big torus $T=\CC^*\otimes N$.
Let $s\colon N\to N'$ be a cosection, that is, $s\circ\iota$ is the identity. Although there is a choice involved here, it will not have any ef\/fect on the end result.

From the above exact sequence, we get two polyhedra in $N_\QQ'$:
\begin{gather*}
\Delta_0 :=s(u^{-1}(1)\cap\sigma),\qquad
\Delta_\infty :=s(u^{-1}(-1)\cap\sigma).
\end{gather*}
Note that both $\Delta_0$ and $\Delta_\infty$ have the tailcone $s( (\ker u) \cap \sigma)$.

We can recover $\cX_0$ from these polyhedra. Let $e$ be a generator of the $\ZZ$ factor in $N' \oplus \ZZ$. Then up to lattice isomorphism, the cone $\sigma$ is equal to
\begin{gather*}
\cone \{\tail(\Delta_0),\Delta_0+e,\Delta_\infty-e\}\subset(N'\oplus \ZZ)_\QQ.
\end{gather*}

We now come to the third choice: let $\Delta_0^0$, $\Delta_0^1$ be an admissible pair of polyhedra such that $\Delta_0^0+\Delta_0^1=\Delta_0$, and such that the pair $\Delta_0^1$, $\Delta_\infty$ is also admissible.

\begin{Theorem}[{cf.~\cite[Theorem~2.8]{ilten:09b}}]\label{thm:family}
Under the above assumptions, there is a flat family \mbox{$\pi\colon \cX \to\PP^1$} with $\cX_0$ the fiber over the origin. The small torus $T'$ acts on $\cX$, preserving fibers and extending the action on $\cX_0$. The fiber $\cX_\infty$ over $\infty$ is toric and isomorphic to $\tv(\sigma_\infty)$, where
\begin{gather*}
\sigma_\infty:=\cone \big\{\tail(\Delta_0),\Delta_0^0+e,\Delta_0^1+\Delta_\infty-e\big\}\subset(N'\oplus \ZZ)_\QQ.
\end{gather*}
Note that up to lattice isomorphism, $\sigma_\infty$ doesn't depend on the cosection $s$.
\end{Theorem}

\begin{proof}
	This is an immediate consequence of the construction of Section~2 in~\cite{ilten:09b}, coupled with the observation that the af\/f\/ine base $B$ of the family there may be replaced by $\PP^1$.
	More precisely, the polyhedra $\Delta_0$ and $\Delta_\infty$ form the coef\/f\/icients of a \emph{$p$-divisor} $\D$ encoding $\cX_0$ as a $T'$-variety, see Remark~1.8 of loc.~cit.  Applying \cite[Theorem~2.8]{ilten:09b} to the decomposition $\Delta_0=\Delta_0^0+\Delta_0^1$ gives a one-parameter f\/lat family over $\mathbb{A}^1$ with $\cX_0$ as the special f\/iber. Since the pair $\Delta_0^1$, $\Delta_\infty$ is admissible, this family may in fact be extended to $\PP^1$, with the f\/iber over $\infty$ described by the $p$-divisor with coef\/f\/icients $\Delta_0^0$ and $\Delta_0^1+\Delta_\infty$. But this is just the toric variety $\tv(\sigma_\infty)$.

	The remark concerning the independence from the cosection $s$ follows from the fact that choosing a dif\/ferent cosection $s'$ will shift $\Delta_0$ by some lattice element $v\in N'$ and $\Delta_\infty$ by its opposite $-v$.
\end{proof}

\begin{Remark}\label{rem:gen}
	If $\Delta_\infty$ is an integral translate of its tailcone, then the general f\/iber of $\pi$ is isomorphic to $\cX_\infty$. However, in general, the general f\/iber will not be toric, admitting instead only the codimension-one torus action by $T'$. The f\/ibers of the family as well as its total space may be described quite explicitly using \emph{$p$-divisors}, see~\cite{altmann:06a}.
\end{Remark}

There are two ways to generalize the above result to non-af\/f\/ine toric varieties. First of all, instead of considering cones and polyhedra, one may consider fans and polyhedral complexes; see \cite[Section~4]{ilten:09b}. A similar approach may also be found in~\cite{mavlyutov:11a}.
We will not discuss this here.
Secondly, if we restrict to projective toric varieties, we can put them in f\/lat projective families by considering degree zero deformations of their af\/f\/ine cones. Explicitly, if $\cX_0=\Proj \CC[\sigma^\vee\cap M]$ with~$\ZZ$ grading given by $v\in N$, this means that we must choose~$u$ as above such that $u(v)=0$. This is the approach which we will pursue in the remainder of this article.
\begin{figure}[t]
\centering
\includegraphics{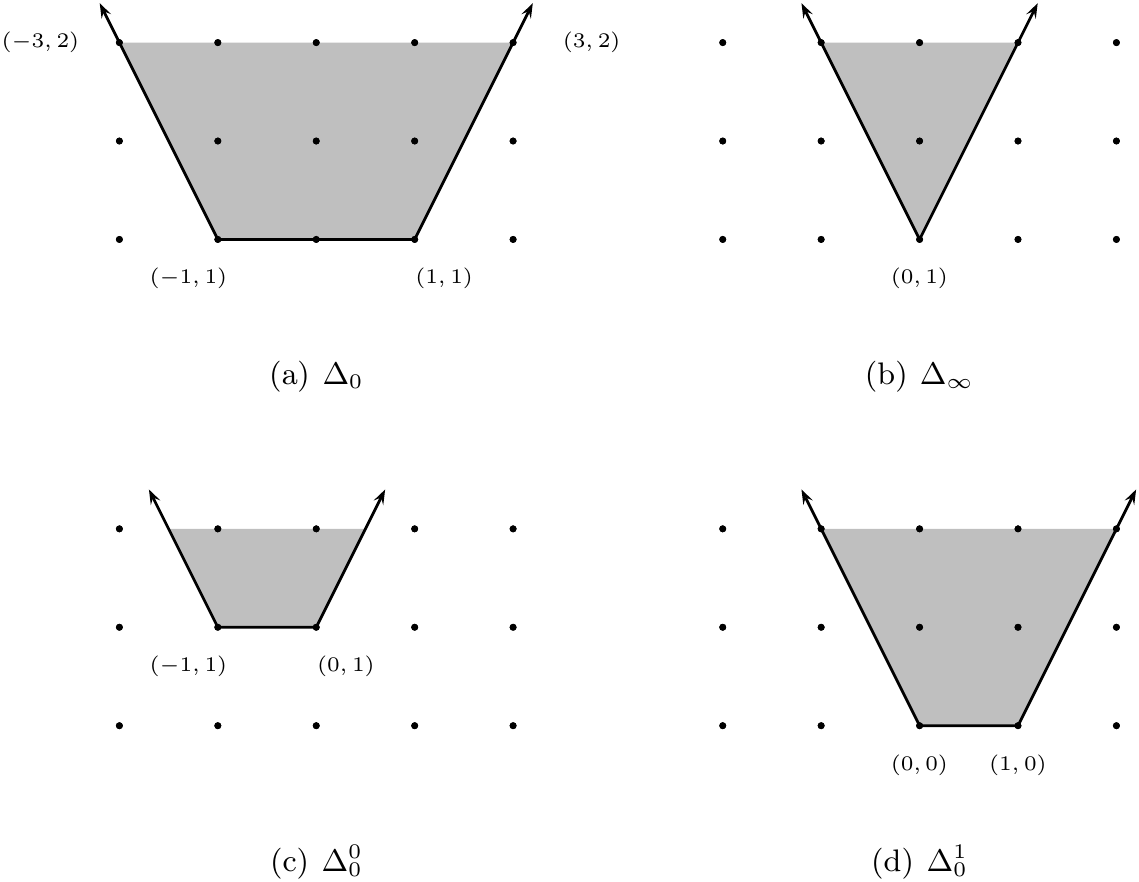}
\caption{Smoothing $\PP(1,1,2)$ to $\PP^1\times\PP^1$.}\label{fig:p112}
\end{figure}

\begin{Example}\label{ex:p112}
Let $N=\ZZ^3$, and consider the cone $\sigma$ generated by $(-1,1,1)$, $(1,1,1)$, and $(0,-1,1)$.
The toric variety $\tv(\sigma)$ is just the af\/f\/ine cone over $\PP(1,1,2)$ in its anticanonical embedding, where the grading is given by $v=(0,0,1)$. Now take $u=(0,1,0)$ as in the construction above. Taking the natural cosection coming from the splitting of $\ZZ^3$, the polyhedra~$\Delta_0$,~$\Delta_\infty$ are as pictured in Fig.~\ref{fig:p112}. There, we also picture polyhedra~$\Delta_0^0$,~$\Delta_0^1$ which are admissible and sum to $\Delta_0$, and such that $\Delta_0^1$, $\Delta_\infty$ are also admissible.
The corresponding cone $\sigma_\infty$ is generated by $(-1,1,1)$, $(0,1,1)$, $(0,-1,1)$, and $(1,-1,1)$. The toric variety $\tv(\sigma)$ is just the af\/f\/ine cone over $\PP^1\times\PP^1$ in its anticanonical embedding, where the grading is again given by $v=(0,0,1)$.

Thus we get a f\/lat family $\cX\to\PP^1$ with $\PP(1,1,2)$ and  $\PP^1\times\PP^1$ as special f\/ibers. Note that since $\Delta_\infty$ is a lattice translate of its tailcone, the general f\/iber of this family is also $\PP^1\times\PP^1$. In fact, this is just the classical smoothing of $\PP(1,1,2)$ to $\PP^1\times\PP^1$.
\end{Example}

\section{Proof of main theorem}\label{sec:main}

In this section, we use Theorem \ref{thm:family} to prove Theorem \ref{thm:main}. Let $g$ be as in Example \ref{ex:easymut} with corresponding birational transformation $\phi$, and suppose that $\phi$ is a mutation of some Laurent polynomial $f$. Then as previously noted,
$f$ can be written as
\begin{gather*}
f=\sum_{i=k}^l f_i z_1^i
\end{gather*}
with $f_i$ Laurent polynomials in $z_2,\ldots,z_n$ such that, for $i>0$, $f_i/g^i$ is a Laurent polynomial. Since we are assuming that $\Delta(f)$ contains the origin in its interior, we must have $k<0<l$.

The coordinates $z_1,\ldots,z_n$ correspond to a basis $e_1,\ldots,e_n$ of $\ZZ^n$.
Let $N=\ZZ\oplus \ZZ^n$ with basis $e_0,\ldots,e_n$, and let $\sigma=\cone\{\Delta(f)+e_0\}$. Then $\tv(\sigma)$ is the af\/f\/ine cone over $\tv(\Delta(f))$, with grading given by $e_0$.

Now, take $u=e_1^*$, and let $s:N\to N'$ be the cosection coming from the induced splitting of $N$.
Let $\Delta_0$ and $\Delta_\infty$ be as in Section~\ref{sec:fam}
and let $\tau= s( (\ker u)\cap \sigma)$ be their common tailcone.
Then $\Delta_0$ is the Minkowski sum of $\tau$ with the  polytope
\begin{gather*}
\conv \left \{\frac{\Delta(f_i)}{i}+\frac{e_0}{i}\right\}_{i=1}^l.
\end{gather*}
Likewise, $\Delta_\infty$ is the Minkowski sum of $\tau$ with the polytope
\begin{gather*}
\conv \left\{\frac{\Delta(f_i)}{-i}+\frac{e_0}{-i} \right\}_{i=k}^{-1}.
\end{gather*}
Taking
\begin{gather*}
\Delta_0^0 :=\tau+\conv \left \{\frac{\Delta(f_i/g^i)}{i}+\frac{e_0}{i}\right\}_{i=1}^l,\qquad
\Delta_0^1 :=\tau+\Delta(g)
\end{gather*}
gives an admissible pair such that $\Delta_0^0+\Delta_0^1=\Delta_0$, and $\Delta_0^1$, $\Delta_\infty$ is admissible as well. Indeed, admissibility follows from the fact that $\Delta_0^1$ is a lattice polyhedron, and the identity $\Delta_0^0+\Delta_0^1=\Delta_0$ follows from the distributive law for convex hulls and Minkowski sums.

We may thus use Theorem~\ref{thm:family} to get a f\/lat family $\pi:\cX\to\PP^1$ with special f\/ibers $\cX_0=\tv(\sigma)$ and $\cX_\infty$.
The cone $\sigma_\infty$ describing $\cX_\infty$ is the cone generated by $\tau$ and by elements of
\begin{gather*}
\Delta\big(f_i/g^i\big)+e_0+ie_1
\end{gather*}
for $k\leq i \leq l$, $i\neq 0$.

Now, note that
\begin{gather*}
\phi(f)=\sum_{i=k}^l \big(f_i/g^i\big)z_1^i.
\end{gather*}
Let $\sigma'=\cone\{\Delta(\phi(f))+e_0\}$. In a moment, we will show that $s( (\ker u)\cap \sigma')=\tau$. It then follows from the description of $\sigma_\infty$ that $\sigma_\infty=\sigma'$. Applying $\Proj$ to the family $\pi$ thus gives us a f\/lat projective family with special f\/ibers $\tv(\Delta(f))$ and $\tv(\Delta(\phi(f)))$, proving Theorem~\ref{thm:main}.

It remains to be shown that $s( (\ker u)\cap \sigma')=\tau$. We will only show the inclusion $\tau \subset s( (\ker u)\cap \sigma')$; the inclusion $s( (\ker u)\cap \sigma')\subset \tau$ follows from an almost identical argument. So let us consider some $v\in\tau$. For $k\leq i \leq l$ we can f\/ind $\lambda_i\in\QQ_{\geq 0}$ and $v_i\in\Delta(f_i)$ satisfying $\sum i\lambda_i=0$ such that
\begin{gather*}
v=\sum_{i=k}^l \lambda_i v_i.
\end{gather*}
For $i>0$, let $v_i'\in\Delta(f_i/g^i)$, $w_i\in\Delta(g^i)$ be such that $v_i=v_i'+w_i$.
Set
\begin{gather*}
w=\frac{1}{\sum\limits_{i=1}^l i \lambda_i}
\sum\limits_{i=1}^l{\lambda_i w_i}\end{gather*}
and note that $w\in\Delta(g)$.
Indeed, for each $i>0$, $\frac{w_i}{i}\in\Delta(g)$, and $w$ is in their convex hull.
Thus
\begin{gather*}
v=\sum_{i=1}^l \lambda_i v_i'+\sum_{i=k}^0 \lambda_i(v_i-iw)
\end{gather*}
which clearly lies in $s( (\ker u)\cap \sigma')$. This completes the proof of Theorem \ref{thm:main}.
\begin{Example}
	Consider the Laurent polynomial $f=x^{-1}y+2y+xy+y^{-1}$. The mutation $x\mapsto x$, $y\mapsto y/(1+x)$ sends $f$ to $f'=x^{-1}y+y+y^{-1}+xy^{-1}$. We have $\tv(\Delta(f))=\PP(1,1,2)$ and $\tv(\Delta(f'))=\PP^1\times\PP^1$, and corresponding family of Theorem \ref{thm:main} is exactly that discussed in Example~\ref{ex:p112}.
\end{Example}

\section{Remarks on mutations and surfaces}\label{sec:surfaces}

We now discuss the connection between Theorem~\ref{thm:main} and the mutations of \cite{galkin}.
Let $N$ be a rank two lattice with dual $M$, and consider some primitive element $u\in M\setminus\{0\}$. Let $e_1$, $e_2$ be a basis of $N$ such that $u=e_2^*$, and let $x$, $y$ be the corresponding monomials in $\CC[N]$, i.e.\ $x=\chi^{e_1}$ and $y=\chi^{e_2}$. Galkin and Usnich consider the birational transformation
\begin{gather*}
\phi_u: \  \CC(N) \to\CC(N), \qquad x\mapsto x, \qquad y\mapsto \frac{y}{1+x}.
\end{gather*}
\begin{Remark}
Dif\/ferent choices of basis will arise in transformations dif\/fering by an element $\psi\in\aut(\CC(N))$ of the form $\psi(\chi^u)=\chi^{A(u)}$, where $A$ is an automorphism of $N$. Such a~map~$\psi$ always maps any Laurent polynomial~$f$ to a Laurent polynomial, and~$\Delta(f)$ and~$\Delta(\psi(f))$ are lattice isomorphic. Thus, we will henceforth ignore the role that the choice of basis plays.
\end{Remark}

\begin{figure}[t]
\centering
\includegraphics{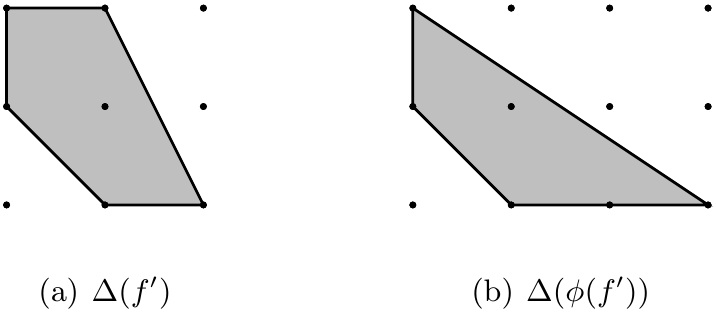}
\caption{A mutation of $f'$.}\label{fig:ex1}
\end{figure}
\begin{figure}[t]
\centering
\includegraphics{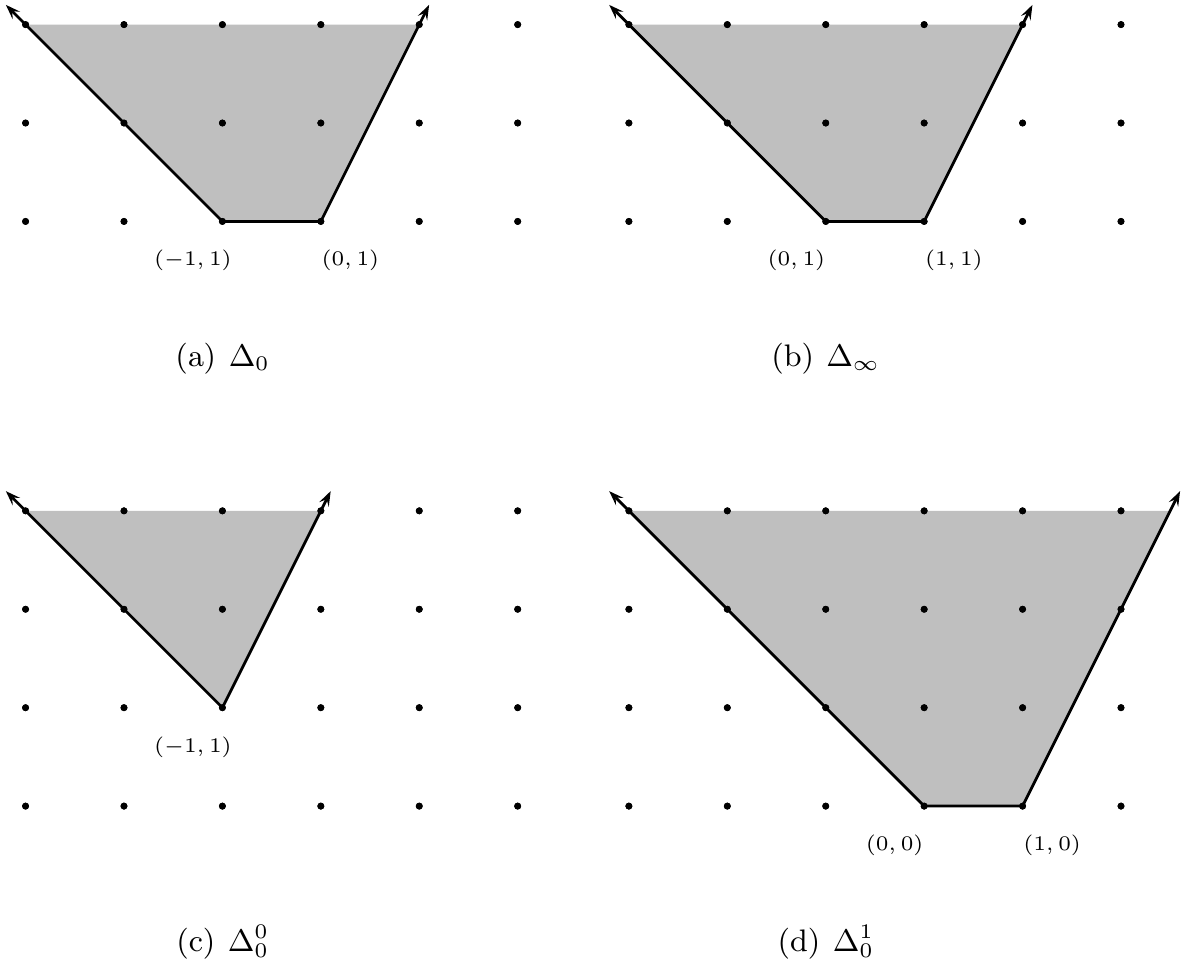}
\caption{A family corresponding to a mutation.}\label{fig:ex2}
\end{figure}

\looseness=-1
Now, let $\Delta$ be any lattice polytope in~$N_\QQ$ containing the origin whose vertices are all pri\-mi\-tive lattice vectors. Such polytopes describe all toric del Pezzo surfaces via~$\tv(\Delta)$. For any facet $\tau\prec \Delta$, let $u_{\tau}$ be the primitive element of~$M$ whose maximum value on~$\Delta$ is achieved exactly on~$\tau$. Note that~$u(\tau)$ measures the lattice height of~$\tau$. Def\/ine $\phi_\tau$ to be the birational transformation~$\phi_{u_\tau}$ from above.

\looseness=-1
In \cite[Section~3]{galkin} there is a list of ten Laurent polynomials whose constant terms series correspond to the Gromov--Witten theory of the ten dif\/ferent families of smooth del Pezzo surfaces. For each polynomial $f$ in the list, $\tv(\Delta(f))$ is a singular toric del Pezzo surface admitting a~$\QQ$-Gorenstein smoothing. Furthermore, $\phi_\tau$ is a mutation for each facet $\tau$, and this property holds for the new Laurent polynomial $\phi_\tau(f)$ and subsequent mutations thereof as well, see \cite[Theorem~16]{galkin}. For each polynomial~$f$, we can thus def\/ine an inf\/inite graph $\Gamma_f$ whose vertices are Laurent polynomials attained from~$f$ via such mutations, and whose edges correspond to mutations.

	To any edge of one of the above graphs $\Gamma_f$, Theorem~\ref{thm:main} says that we can associate a~deformation between the toric varieties corresponding to that edge's vertices.

\begin{Example}
	Consider the Laurent polynomial $f=x+y+x^{-1}y^{-1}+x^{-1}+y^{-1}$ from the list of \cite{galkin} corresponding to the del Pezzo surface of degree seven.
	This is lattice equivalent to the polynomial $f'=x^{-1}+x^{-1}y+y+y^{-1}+xy^{-1}$. Mutating this polynomial (with respect to this choice of~$x$,~$y$) gives $\phi(f')=x^{-1}+x^{-1}y+y^{-1}+xy^{-1}+x^2y^{-1}$. In Fig.~\ref{fig:ex1}, we picture the Newton polytopes of $f'$ and $\phi(f')$. In Fig.~\ref{fig:ex2} we picture the polytopes $\Delta_0$, $\Delta_\infty$, $\Delta_0^0$, and $\Delta_0^1$ as in the proof of Theorem~\ref{thm:main}.
\end{Example}

\pdfbookmark[1]{References}{ref}
\LastPageEnding

\end{document}